\documentclass[a4paper, 12 pt]{article}

\usepackage[T2A]{fontenc}
\usepackage[cp1251]{inputenc}
\usepackage[english]{babel}
\usepackage[tbtags]{amsmath}
\usepackage{amsfonts,amssymb}

\usepackage{mathrsfs}

\usepackage{geometry} 
\begin{document}
\title{Fractal generalized Pascal matrices}
 
 \author{E. Burlachenko}
 \date{}
 
 \maketitle

\begin{abstract}

Set of generalized Pascal matrices whose elements are generalized binomial coefficients is considered as an integral object. The special system of generalized Pascal matrices, based on which we are building fractal generalized Pascal matrices,  is introduced. Pascal matrix (Pascal triangle) is the Hadamard product of the fractal generalized Pascal matrices whose elements equal to ${{p}^{k}}$, where $p$ is a fixed prime number, $k=0$, $1$, $2$, … .

The concept of zero generalized Pascal matrices, an example of which is the Pascal triangle modulo 2, arise  in connection with the system of matrices introduced.

\end{abstract}
\section{Introduction}

Consider the following generalization of the binomial coefficients [1]. For the coefficients of the formal power series  $b\left( x \right)$, ${{b}_{0}}=0$; ${{b}_{n}}\ne 0$,  $n>0$, denote
$${{b}_{0}}!=1,\qquad {{b}_{n}}!=\prod\limits_{m=1}^{n}{{{b}_{m}}},\qquad {{\left( \begin{matrix}
   n  \\
   m  \\
\end{matrix} \right)}_{b}}=\frac{{{b}_{n}}!}{{{b}_{m}}!{{b}_{n-m}}!};\qquad   {{\left( \begin{matrix}
   n  \\
   m  \\
\end{matrix} \right)}_{b}}=0,  \qquad m>n.$$
Then
\[{{\left( \begin{matrix}
   n  \\
   m  \\
\end{matrix} \right)}_{b}}={{\left( \begin{matrix}
   n-1  \\
   m-1  \\
\end{matrix} \right)}_{b}}+\frac{{{b}_{n}}-{{b}_{m}}}{{{b}_{n-m}}}{{\left( \begin{matrix}
   n-1  \\
   m  \\
\end{matrix} \right)}_{b}}.\]
$n$-th coefficient of the series $a\left( x \right)$, $\left( n,m \right)$-th element of the matrix $A$, $n$-th row and $n$-th column of the matrix$A$ will be denoted  respectively by 
$$\left[ {{x}^{n}} \right]a\left( x \right), \qquad {{\left( A \right)}_{n,m}}, \qquad  \left[ n,\to  \right]A,  \qquad \left[\uparrow ,n \right]A.$$
We associate rows and columns of matrices with the generating functions of their elements. For the elements of the lower triangular matrices will be appreciated that ${{\left( A \right)}_{n,m}}=0$, if $n<m$. 
Consider matrix
$${{P}_{c\left( x \right)}}=\left( \begin{matrix}
   \frac{{{c}_{0}}{{c}_{0}}}{{{c}_{0}}} & 0 & 0 & 0 & 0 & \ldots   \\
   \frac{{{c}_{0}}{{c}_{1}}}{{{c}_{1}}} & \frac{{{c}_{1}}{{c}_{0}}}{{{c}_{1}}} & 0 & 0 & 0 & \ldots   \\
   \frac{{{c}_{0}}{{c}_{2}}}{{{c}_{2}}} & \frac{{{c}_{1}}{{c}_{1}}}{{{c}_{2}}} & \frac{{{c}_{2}}{{c}_{0}}}{{{c}_{2}}} & 0 & 0 & \ldots   \\
   \frac{{{c}_{0}}{{c}_{3}}}{{{c}_{3}}} & \frac{{{c}_{1}}{{c}_{2}}}{{{c}_{3}}} & \frac{{{c}_{2}}{{c}_{1}}}{{{c}_{3}}} & \frac{{{c}_{3}}{{c}_{0}}}{{{c}_{3}}} & 0 & \ldots   \\
   \frac{{{c}_{0}}{{c}_{4}}}{{{c}_{4}}} & \frac{{{c}_{1}}{{c}_{3}}}{{{c}_{4}}} & \frac{{{c}_{2}}{{c}_{2}}}{{{c}_{4}}} & \frac{{{c}_{3}}{{c}_{1}}}{{{c}_{4}}} & \frac{{{c}_{4}}{{c}_{0}}}{{{c}_{4}}} & \ldots   \\
   \vdots  & \vdots  & \vdots  & \vdots  & \vdots  & \ddots   \\
\end{matrix} \right).$$
$${{\left( {{P}_{c\left( x \right)}} \right)}_{n,m}}=\frac{{{c}_{m}}{{c}_{n-m}}}{{{c}_{n}}},  \qquad{{c}_{n}}\in \mathbb{R}, \qquad{{c}_{n}}\ne 0.$$
 Denote $[\uparrow ,1]{{P}_{c\left( x \right)}}=b\left( x \right)$. If ${{c}_{0}}=1$,  then

$${{c}_{n}}=\frac{c_{1}^{n}}{{{b}_{n}}!} ,\qquad {{\left( {{P}_{c\left( x \right)}} \right)}_{n,m}}={{\left( \begin{matrix}
   n  \\
   m  \\
\end{matrix} \right)}_{b}}.$$
Let ${{c}_{0}}=1$. Since ${{P}_{c\left( x \right)}}={{P}_{c\left( \varphi x \right)}}$, we take for uniqueness that ${{c}_{1}}=1$. Matrix ${{P}_{c\left( x \right)}}$ will be called generalized Pascal matrix. If $c\left( x \right)={{e}^{x}}$, it turns into Pascal matrix
$$P=\left( \begin{matrix}
   1 & 0 & 0 & 0 & 0 & \ldots   \\
   1 & 1 & 0 & 0 & 0 & \ldots   \\
   1 & 2 & 1 & 0 & 0 & \ldots   \\
   1 & 3 & 3 & 1 & 0 & \ldots   \\
   1 & 4 & 6 & 4 & 1 & \ldots   \\
   \vdots  & \vdots  & \vdots  & \vdots  & \vdots  & \ddots   \\
\end{matrix} \right).$$

In Sec. 2. we consider the set of generalized Pascal matrices as a group under  Hadamard multiplication and introduce a special system of matrices, which implies the concept of zero generalized Pascal matrices; as an example, we consider the matrix of the $q$-binomial coefficients for $q=-1$. In Sec. 3. we construct a fractal generalized Pascal matrices whose Hadamard product is the Pascal matrix. In Sec. 4. we consider the fractal zero generalized Pascal matrices that were previously considered in [7] – [11] (generalized  Sierpinski matrices and others).

\section{Special system of generalized Pascal matrices}

Elements of the matrix ${{P}_{c\left( x \right)}}$, – denote them ${{\left( {{P}_{c\left( x \right)}} \right)}_{n,m}}={{{c}_{m}}{{c}_{n-m}}}/{{{c}_{n}}}\;=\left( n,m \right)$ for generality which will be discussed later, –  satisfy the identities 
	$$\left( n,0 \right)=1, \qquad\left( n,m \right)=\left( n,n-m \right), \eqno  	(1)$$
$$\text{  }\left( n+q,q \right)\left( n+p,m+p \right)\left( m+p,p \right)=\left( n+p,p \right)\left( n+q,m+q \right)\left( m+q,q \right),\eqno(2)$$
$q$, $p=0$, $1$, $2$, … It means that each matrix ${{P}_{c\left( x \right)}}$ can be associated with the algebra of formal power series whose elements are multiplied by the rule
$$a\left( x \right)\circ b\left( x \right)=g\left( x \right),  \qquad{{g}_{n}}=\sum\limits_{m=0}^{n}{\left( n,m \right){{a}_{m}}}{{b}_{m-n}},$$
that is, if ${{\left( A \right)}_{n,m}}={{a}_{n-m}}\left( n,m \right)$, ${{\left( B \right)}_{n,m}}={{b}_{n-m}}\left( n,m \right)$, ${{\left( G \right)}_{n,m}}={{g}_{n-m}}\left( n,m \right)$,   then $AB=BA=G$:
$${{g}_{n}}={{\left( n+p,p \right)}^{-1}}\sum\limits_{m=0}^{n}{\left( n+p,m+p \right)\left( m+p,p \right){{a}_{n-m}}{{b}_{m}}}=$$
$$={{\left( n+q,q \right)}^{-1}}\sum\limits_{m=0}^{n}{\left( n+q,m+q \right)\left( m+q,q \right){{a}_{n-m}}{{b}_{m}}}.$$

The set of generalized Pascal matrix is a group under Hadamard multiplication (we denote this operation $\times $):
$${{P}_{c\left( x \right)}}\times {{P}_{g\left( x \right)}}={{P}_{c\left( x \right)\times g\left( x \right)}},  \qquad c\left( x \right)\times g\left( x \right)=\sum\limits_{n=0}^{\infty }{{{c}_{n}}{{g}_{n}}}{{x}^{n}}.$$

Introduce the special system of matrices
$${}_{\varphi ,q}P={}_{q}P\left( \varphi  \right)={{P}_{c\left( \varphi ,q,x \right)}},  \qquad c\left( \varphi ,q,x \right)=\left( \sum\limits_{n=0}^{q-1}{{{x}^{n}}} \right){{\left( 1-\frac{{{x}^{q}}}{\varphi } \right)}^{-1}}, \qquad q>1.$$
Then 
$${{c}_{qn+i}}=\frac{1}{{{\varphi }^{n}}},\qquad  0\le i<q; \qquad {{c}_{qn-i}}=\frac{1}{{{\varphi }^{n-1}}},\qquad  0<i\le q,$$
$$\frac{{{c}_{qm+j}}{{c}_{q\left( n-m \right)+i-j}}}{{{c}_{qn+i}}}=\frac{{{\varphi }^{n}}}{{{\varphi }^{m}}{{\varphi }^{n-m}}}=1,\qquad i\ge j;  \qquad=\frac{{{\varphi }^{n}}}{{{\varphi }^{m}}{{\varphi }^{n-m-1}}}=\varphi , \qquad i<j,$$
or
$${{\left( _{\varphi ,q}P \right)}_{n,m}}=1, \qquad n\left( \bmod q \right)\ge m\left( \bmod q \right); \qquad=\varphi , \qquad n\left( \bmod q \right)<m\left( \bmod q \right).$$
For example, $_{\varphi ,2}P$, $_{\varphi ,3}P$:
$$\left( \begin{matrix}
   1 & 0 & 0 & 0 & 0 & 0 & 0 & 0 & 0 & \ldots   \\
   1 & 1 & 0 & 0 & 0 & 0 & 0 & 0 & 0 & \ldots   \\
   1 & \varphi  & 1 & 0 & 0 & 0 & 0 & 0 & 0 & \ldots   \\
   1 & 1 & 1 & 1 & 0 & 0 & 0 & 0 & 0 & \ldots   \\
   1 & \varphi  & 1 & \varphi  & 1 & 0 & 0 & 0 & 0 & \ldots   \\
   1 & 1 & 1 & 1 & 1 & 1 & 0 & 0 & 0 & \dots   \\
   1 & \varphi  & 1 & \varphi  & 1 & \varphi  & 1 & 0 & 0 & \ldots   \\
   1 & 1 & 1 & 1 & 1 & 1 & 1 & 1 & 0 & \ldots   \\
   1 & \varphi  & 1 & \varphi  & 1 & \varphi  & 1 & \varphi  & 1 & \ldots   \\
   \vdots  & \vdots  & \vdots  & \vdots  & \vdots  & \vdots  & \vdots  & \vdots  & \vdots  & \ddots   \\
\end{matrix} \right), \qquad \left( \begin{matrix}
   1 & 0 & 0 & 0 & 0 & 0 & 0 & 0 & 0 & \ldots   \\
   1 & 1 & 0 & 0 & 0 & 0 & 0 & 0 & 0 & \ldots   \\
   1 & 1 & 1 & 0 & 0 & 0 & 0 & 0 & 0 & \ldots   \\
   1 & \varphi  & \varphi  & 1 & 0 & 0 & 0 & 0 & 0 & \ldots   \\
   1 & 1 & \varphi  & 1 & 1 & 0 & 0 & 0 & 0 & \ldots   \\
   1 & 1 & 1 & 1 & 1 & 1 & 0 & 0 & 0 & \ldots   \\
   1 & \varphi  & \varphi  & 1 & \varphi  & \varphi  & 1 & 0 & 0 & \ldots   \\
   1 & 1 & \varphi  & 1 & 1 & \varphi  & 1 & 1 & 0 & \ldots   \\
   1 & 1 & 1 & 1 & 1 & 1 & 1 & 1 & 1 & \ldots   \\
   \vdots  & \vdots  & \vdots  & \vdots  & \vdots  & \vdots  & \vdots  & \vdots  & \vdots  & \ddots   \\
\end{matrix} \right).$$

Elements of the matrix $_{\varphi ,q}P\times {{P}_{c\left( x \right)}}$ satisfy the identities (1), (2) for any values $\varphi $, so it makes sense to consider also the case $\varphi =0$ since it corresponds to a certain algebra of formal power series (which, obviously, contains zero divisors; for example, in the algebra associated with the matrix $_{0,2}P\times {{P}_{c\left( x \right)}}$ the product of the series of the form $xa\left( {{x}^{2}} \right)$ is zero). It is clear that in this case the series $c\left( \varphi ,q,x \right)$ is not defined. Matrix $_{0,q}P\times {{P}_{c\left( x \right)}}$ and Hadamard product of such matrices will be called zero generalized Pascal matrix.

{\bfseries Remark.} Zero generalized Pascal matrix appears when considering the set of generalized Pascal matrices ${{P}_{g\left( q,x \right)}}$:
$${{\left( {{P}_{g\left( q,x \right)}} \right)}_{n,1}}=\left[ {{x}^{n}} \right]\frac{x}{\left( 1-x \right)\left( 1-qx \right)}=\sum\limits_{m=0}^{n-1}{{{q}^{m}}}, \qquad q\in \mathbb{R}.$$
Here $g\left( 0,x \right)={{\left( 1-x \right)}^{-1}}$,  $g\left( 1,x \right)={{e}^{x}}$. In other cases  (the $q$-umbral calculus [2]), except $q=-1$,
$$g\left( q,x \right)=\sum\limits_{n=0}^{\infty }{\frac{{{\left( q-1 \right)}^{n}}}{\left( {{q}^{n}}-1 \right)!}}{{x}^{n}},  \qquad\left( {{q}^{n}}-1 \right)!=\prod\limits_{m=1}^{n}{\left( {{q}^{m}}-1 \right)},  \qquad\left( {{q}^{0}}-1 \right)!=1.$$
Matrices  ${{P}_{g\left( q,x \right)}}$, $P_{g\left( q,x \right)}^{-1}$ also can  be defined as follows:
$$[\uparrow ,n]{{P}_{g\left( q,x \right)}}={{x}^{n}}\prod\limits_{m=0}^{n}{{{\left( 1-{{q}^{m}}x \right)}^{-1}}},    \qquad[n,\to ]P_{g\left( q,x \right)}^{-1}=\prod\limits_{m=0}^{n-1}{\left( x-{{q}^{m}} \right)}.$$
When  $q=-1$ we get the matrices ${{P}_{g\left( -1,x \right)}}$, $P_{g\left( -1,x \right)}^{-1}$:
$$\left( \begin{matrix}
   1 & 0 & 0 & 0 & 0 & 0 & 0 & \ldots   \\
   1 & 1 & 0 & 0 & 0 & 0 & 0 & \ldots   \\
   1 & 0 & 1 & 0 & 0 & 0 & 0 & \ldots   \\
   1 & 1 & 1 & 1 & 0 & 0 & 0 & \ldots   \\
   1 & 0 & 2 & 0 & 1 & 0 & 0 & \ldots   \\
   1 & 1 & 2 & 2 & 1 & 1 & 0 & \ldots   \\
   1 & 0 & 3 & 0 & 3 & 0 & 1 & \ldots   \\
   \vdots  & \vdots  & \vdots  & \vdots  & \vdots  & \vdots  & \vdots  & \ddots   \\
\end{matrix} \right), \qquad\left( \begin{matrix}
   1 & 0 & 0 & 0 & 0 & 0 & 0 & \ldots   \\
   -1 & 1 & 0 & 0 & 0 & 0 & 0 & \ldots   \\
   -1 & 0 & 1 & 0 & 0 & 0 & 0 & \ldots   \\
   1 & -1 & -1 & 1 & 0 & 0 & 0 & \ldots   \\
   1 & 0 & -2 & 0 & 1 & 0 & 0 & \ldots   \\
   -1 & 1 & 2 & -2 & -1 & 1 & 0 & \ldots   \\
   -1 & 0 & 3 & 0 & -3 & 0 & 1 & \ldots   \\
   \vdots  & \vdots  & \vdots  & \vdots  & \vdots  & \vdots  & \vdots  & \ddots   \\
\end{matrix} \right),$$
where the series $g\left( -1,x \right)$  is not defined. Since
$${{\left( {{P}_{g\left( -1,x \right)}} \right)}_{2n+i,2m+j}}=\left[ {{x}^{2n+i}} \right]\frac{{{\left( 1+x \right)}^{1-j}}{{x}^{2m+j}}}{{{\left( 1-{{x}^{2}} \right)}^{m+1}}}=\left( \begin{matrix}
   n  \\
   m  \\
\end{matrix} \right), \qquad i\ge j; \qquad=0, \qquad i<j;$$ 
$i,j=0,$$1$, then
$${{P}_{g\left( -1,x \right)}}={}_{0,2}P\times {{P}_{c\left( x \right)}},  \qquad c\left( x \right)=\left( 1+x \right){{e}^{{{x}^{2}}}}:$$
$${{c}_{2n+i}}=\frac{1}{n!},  \qquad 0\le i<2;  \qquad{{c}_{2n-i}}=\frac{1}{\left( n-1 \right)!},  \qquad 0<i\le 2,$$
$${{\left( {{P}_{c\left( x \right)}} \right)}_{2n+i,2m+j}}=\left( \begin{matrix}
   n  \\
   m  \\
\end{matrix} \right), \qquad i\ge j; \qquad=n\left( \begin{matrix}
   n-1  \\
   m  \\
\end{matrix} \right), \qquad i<j.$$
A generalization of this matrix is the matrix 
$$_{0,q}P\times {{P}_{c\left( x \right)}}, \qquad c\left( x \right)=\left( \sum\limits_{n=0}^{q-1}{{{x}^{n}}} \right){{e}^{{{x}^{q}}}},$$
$${{\left( {{P}_{c\left( x \right)}} \right)}_{qn+i,qm+j}}=\left( \begin{matrix}
   n  \\
   m  \\
\end{matrix} \right), \qquad i\ge j; \qquad=n\left( \begin{matrix}
   n-1  \\
   m  \\
\end{matrix} \right), \qquad i<j; \qquad 0\le i,j<q.$$

Each nonzero generalized Pascal matrix is the Hadamard product of the matrices $_{\varphi ,q}P$. Since the first column of the matrix ${{P}_{c\left( x \right)}}$, –  denote it $b\left( x \right)$, –  is the Hadamard product of  the first columns of the  matrices $_{\varphi ,q}P$, –  denote them $_{\varphi ,q}b\left( x \right)$:
$$\left[ {{x}^{n}} \right]{}_{\varphi ,q}b\left( x \right)=1, \qquad n\left( \bmod q \right)\ne 0; \qquad=\varphi , \qquad n\left( \bmod q \right)=0,$$
so that
$${{P}_{c\left( x \right)}}={}_{2}P\left( {{b}_{2}} \right)\times {}_{3}P\left( {{b}_{3}} \right)\times {}_{4}P\left( {{{b}_{4}}}/{{{b}_{2}}}\; \right)\times {}_{5}P\left( {{b}_{5}} \right)\times {}_{6}P\left( {{{b}_{6}}}/{{{b}_{2}}{{b}_{3}}}\; \right)\times {}_{7}P\left( {{b}_{7}} \right)\times $$
$${{\times }_{8}}P\left( {{{b}_{8}}}/{{{b}_{4}}}\; \right)\times {}_{9}P\left( {{{b}_{9}}}/{{{b}_{3}}}\; \right)\times {}_{10}P\left( {{{b}_{10}}}/{{{b}_{2}}{{b}_{5}}}\; \right)\times {}_{11}P\left( {{b}_{11}} \right)\times {}_{12}P\left( {{{b}_{12}}{{b}_{2}}}/{{{b}_{4}}{{b}_{6}}}\; \right)\times ...$$ 
and so on. Let ${{e}_{q}}$ is a basis vector of an infinite-dimensional vector space. Mapping of the set of generalized Pascal matrices in an infinite-dimensional vector space such that $_{\varphi ,q}P\to {{e}_{q}}\log \left| \varphi  \right|$ is a group  homomorphism whose kernel consists of all involutions in the group of generalized Pascal matrices, i.e. from matrices whose non-zero elements equal to $\pm 1$. Thus, the set of generalized Pascal matrices whose elements are non-negative numbers is an infinite-dimensional vector space. Zero generalized Pascal matrices can be viewed as points at infinity of space.

\section{Fractal generalized Pascal matrices}

Matrices, which will be discussed (precisely, isomorphic to them), are considered in [3] ($p$-index Pascal triangle), [4], [5, p. 80-88]. These matrices are introduced explicitly in [6] in connection with the generalization of the theorems on the divisibility of binomial coefficients. We consider them from point of view based on the system of matrices $_{\varphi ,q}P$. We start with the matrix
$$_{\left[ 2 \right]}P={}_{2,2}P\times {}_{{{2,2}^{2}}}P\times {}_{{{2,2}^{3}}}P\times ...\times {}_{{{2,2}^{k}}}P\times ...$$
$${}_{\left[ 2 \right]}P=\left(\setcounter{MaxMatrixCols}{20} \begin{matrix}
   1 & 0 & 0 & 0 & 0 & 0 & 0 & 0 & 0 & 0 & 0 & 0 & 0 & 0 & 0 & 0 & \ldots   \\
   1 & 1 & 0 & 0 & 0 & 0 & 0 & 0 & 0 & 0 & 0 & 0 & 0 & 0 & 0 & 0 & \ldots   \\
   1 & 2 & 1 & 0 & 0 & 0 & 0 & 0 & 0 & 0 & 0 & 0 & 0 & 0 & 0 & 0 & \ldots   \\
   1 & 1 & 1 & 1 & 0 & 0 & 0 & 0 & 0 & 0 & 0 & 0 & 0 & 0 & 0 & 0 & \ldots   \\
   1 & 4 & 2 & 4 & 1 & 0 & 0 & 0 & 0 & 0 & 0 & 0 & 0 & 0 & 0 & 0 & \ldots   \\
   1 & 1 & 2 & 2 & 1 & 1 & 0 & 0 & 0 & 0 & 0 & 0 & 0 & 0 & 0 & 0 & \ldots   \\
   1 & 2 & 1 & 4 & 1 & 2 & 1 & 0 & 0 & 0 & 0 & 0 & 0 & 0 & 0 & 0 & \ldots   \\
   1 & 1 & 1 & 1 & 1 & 1 & 1 & 1 & 0 & 0 & 0 & 0 & 0 & 0 & 0 & 0 & \ldots   \\
   1 & 8 & 4 & 8 & 2 & 8 & 4 & 8 & 1 & 0 & 0 & 0 & 0 & 0 & 0 & 0 & \ldots   \\
   1 & 1 & 4 & 4 & 2 & 2 & 4 & 4 & 1 & 1 & 0 & 0 & 0 & 0 & 0 & 0 & \ldots   \\
   1 & 2 & 1 & 8 & 2 & 4 & 2 & 8 & 1 & 2 & 1 & 0 & 0 & 0 & 0 & 0 & \ldots   \\
   1 & 1 & 1 & 1 & 2 & 2 & 2 & 2 & 1 & 1 & 1 & 1 & 0 & 0 & 0 & 0 & \ldots   \\
   1 & 4 & 2 & 4 & 1 & 8 & 4 & 8 & 1 & 4 & 2 & 4 & 1 & 0 & 0 & 0 & \ldots   \\
   1 & 1 & 2 & 2 & 1 & 1 & 4 & 4 & 1 & 1 & 2 & 2 & 1 & 1 & 0 & 0 & \ldots   \\
   1 & 2 & 1 & 4 & 1 & 2 & 1 & 8 & 1 & 2 & 1 & 4 & 1 & 2 & 1 & 0 & \ldots   \\
   1 & 1 & 1 & 1 & 1 & 1 & 1 & 1 & 1 & 1 & 1 & 1 & 1 & 1 & 1 & 1 & \ldots   \\
   \vdots  & \vdots  & \vdots  & \vdots  & \vdots  & \vdots  & \vdots  & \vdots  & \vdots  & \vdots  & \vdots  & \vdots  & \vdots  & \vdots  & \vdots  & \vdots  & \ddots   \\
\end{matrix} \right).$$

The first column of the matrix  $_{\left[ 2 \right]}P$, – denote it $ b\left( x \right)$, – is the Hadamard product of  the series $_{{{2,2}^{k}}}b\left( x \right)$ , $k=1$, $2$ , … ,
$$\left[ {{x}^{n}} \right]{}_{{{2,2}^{k}}}b\left( x \right)=1, \qquad n\left( \bmod {{2}^{k}} \right)\ne 0; \qquad=2, \qquad n\left( \bmod {{2}^{k}} \right)=0.$$
It is the generating function of the distribution of divisors ${{2}^{k}}$ in the series of natural numbers:
$$b\left( x \right)=x+2{{x}^{2}}+{{x}^{3}}+4{{x}^{4}}+{{x}^{5}}+2{{x}^{6}}+{{x}^{7}}+8{{x}^{8}}+{{x}^{9}}+2{{x}^{10}}+{{x}^{11}}+4{{x}^{12}}+$$
$$+{{x}^{13}}+2{{x}^{14}}+{{x}^{15}}+16{{x}^{16}}+{{x}^{17}}+2{{x}^{18}}+{{x}^{19}}+4{{x}^{20}}+{{x}^{21}}+2{{x}^{22}}+{{x}^{23}}+...$$
{\bfseries Theorem.}
$${{b}_{{{2}^{k}}n+i}}={{b}_{i}}, \qquad 0<i<{{2}^{k}}.$$
{\bfseries Proof.}  It's obvious that ${{b}_{2n}}=2{{b}_{n}}$,  ${{b}_{2n+1}}=1$. Then
$$b\left( x \right)=\frac{x}{1-{{x}^{2}}}+2b\left( {{x}^{2}} \right)=\sum\limits_{n=0}^{k-1}{\frac{{{2}^{n}}{{x}^{{{2}^{n}}}}}{1-{{x}^{{{2}^{n+1}}}}}}+{{2}^{k}}b\left( {{x}^{{{2}^{k}}}} \right)=\sum\limits_{n=0}^{\infty }{\frac{{{2}^{n}}{{x}^{{{2}^{n}}}}}{1-{{x}^{{{2}^{n+1}}}}}},$$
or  
$$b\left( x \right)=\sum\limits_{n=0}^{\infty }{_{\left( n \right)}a\left( x \right)}, \qquad_{\left( n \right)}a\left( x \right)=\sum\limits_{m=0}^{\infty }{{{2}^{n}}{{x}^{{{2}^{n}}\left( 2m+1 \right)}}},$$ 
$$_{\left( n \right)}{{a}_{{{2}^{n}}\left( 2m+1 \right)}}={}_{\left( n \right)}{{a}_{{{2}^{k}}p+{{2}^{n}}\left( 2m+1 \right)}}={{2}^{n}}, \qquad k>n,$$
and $_{\left( n \right)}{{a}_{m}}=0$ in other cases. Hence, $_{\left( n \right)}{{a}_{{{2}^{k}}p+i}}={}_{\left( n \right)}{{a}_{i}}$, $0<i<{{2}^{k}}$.

Note identities 
$${{b}_{2n+1}}!={{b}_{2n}}!,  \qquad{{b}_{2n-1}}!={{b}_{2\left( n-1 \right)}}!,$$
$${{b}_{2n}}!={{b}_{2n}}{{b}_{2n-1}}!=2{{b}_{n}}{{b}_{2\left( n-1 \right)}}!=2{{b}_{n}}{{b}_{2\left( n-1 \right)}}{{b}_{2\left( n-1 \right)-1}}!={{2}^{2}}{{b}_{n}}{{b}_{\left( n-1 \right)}}{{b}_{2\left( n-2 \right)}}!=...={{2}^{n}}{{b}_{n}}!,$$
$${{b}_{{{2}^{k}}n}}!={{2}^{{{2}^{k-1}}n}}{{b}_{{{2}^{k-1}}n}}!={{2}^{{{2}^{k-1}}n}}{{2}^{{{2}^{k-2}}n}}{{b}_{{{2}^{k-2}}n}}!=...={{2}^{\left( {{2}^{k}}-1 \right)n}}{{b}_{n}}!.$$

The series $\left( x \right)=\sum\nolimits_{n=0}^{\infty }{{{{x}^{n}}}/{{{b}_{n}}!}\;}$  also is the fractal:
$$c\left( x \right)=1+x+\frac{{{x}^{2}}}{2}+\frac{{{x}^{3}}}{2}+\frac{{{x}^{4}}}{{{2}^{3}}}+\frac{{{x}^{5}}}{{{2}^{3}}}+\frac{{{x}^{6}}}{{{2}^{4}}}+\frac{{{x}^{7}}}{{{2}^{4}}}+\frac{{{x}^{8}}}{{{2}^{7}}}+\frac{{{x}^{9}}}{{{2}^{7}}}+\frac{{{x}^{10}}}{{{2}^{8}}}+\frac{{{x}^{11}}}{{{2}^{8}}}+$$
$$+\frac{{{x}^{12}}}{{{2}^{10}}}+\frac{{{x}^{13}}}{{{2}^{10}}}+\frac{{{x}^{14}}}{{{2}^{11}}}+\frac{{{x}^{15}}}{{{2}^{11}}}+\frac{{{x}^{16}}}{{{2}^{15}}}+\frac{{{x}^{17}}}{{{2}^{15}}}+\frac{{{x}^{18}}}{{{2}^{16}}}+\frac{{{x}^{19}}}{{{2}^{16}}}+\frac{{{x}^{20}}}{{{2}^{18}}}+\frac{{{x}^{21}}}{{{2}^{18}}}+...,$$
$${{c}_{2n}}={{c}_{2n+1}}=\frac{1}{{{2}^{n}}{{b}_{n}}!}=\frac{{{c}_{n}}}{{{2}^{n}}},  \qquad c\left( x \right)=\left( 1+x \right)c\left( \frac{{{x}^{2}}}{2} \right)=\prod\limits_{n=0}^{\infty }{\left( 1+\frac{{{x}^{{{2}^{n}}}}}{{{2}^{{{2}^{n}}-1}}} \right)}.$$

Denote
$${{\left( _{\left[ 2 \right]}P \right)}_{n,m}}=\frac{{{b}_{n}}!}{{{b}_{m}}!{{b}_{n-m}}!}={{\left( \begin{matrix}
   n  \\
   m  \\
\end{matrix} \right)}_{2}}.$$
Since ${{b}_{{{2}^{k}}n+i}}!={{b}_{{{2}^{k}}n}}!{{b}_{i}}!={{2}^{\left( {{2}^{k}}-1 \right)n}}{{b}_{n}}!{{b}_{i}}!$, $0\le i<{{2}^{k}}$, then
$${{\left( \begin{matrix}
   {{2}^{k}}n+i  \\
   {{2}^{k}}m+j  \\
\end{matrix} \right)}_{2}}={{\left( \begin{matrix}
   n  \\
   m  \\
\end{matrix} \right)}_{2}}{{\left( \begin{matrix}
   i  \\
   j  \\
\end{matrix} \right)}_{2}}, \qquad i\ge j.$$
For example, matrices $_{\left[ 2 \right]}P\times {}_{0,2}P$, $_{\left[ 2 \right]}P\times {}_{0,4}P$ have the form
$$\left( \begin{matrix}
   1 & 0 & 0 & 0 & 0 & 0 & 0 & 0 & 0 & 0 & 0 & 0 & \ldots   \\
   1 & 1 & 0 & 0 & 0 & 0 & 0 & 0 & 0 & 0 & 0 & 0 & \ldots   \\
   1 & 0 & 1 & 0 & 0 & 0 & 0 & 0 & 0 & 0 & 0 & 0 & \ldots   \\
   1 & 1 & 1 & 1 & 0 & 0 & 0 & 0 & 0 & 0 & 0 & 0 & \ldots   \\
   1 & 0 & 2 & 0 & 1 & 0 & 0 & 0 & 0 & 0 & 0 & 0 & \ldots   \\
   1 & 1 & 2 & 2 & 1 & 1 & 0 & 0 & 0 & 0 & 0 & 0 & \ldots   \\
   1 & 0 & 1 & 0 & 1 & 0 & 1 & 0 & 0 & 0 & 0 & 0 & \ldots   \\
   1 & 1 & 1 & 1 & 1 & 1 & 1 & 1 & 0 & 0 & 0 & 0 & \ldots   \\
   1 & 0 & 4 & 0 & 2 & 0 & 4 & 0 & 1 & 0 & 0 & 0 & \ldots   \\
   1 & 1 & 4 & 4 & 2 & 2 & 4 & 4 & 1 & 1 & 0 & 0 & \ldots   \\
   1 & 0 & 1 & 0 & 2 & 0 & 2 & 0 & 1 & 0 & 1 & 0 & \ldots   \\
   1 & 1 & 1 & 1 & 2 & 2 & 2 & 2 & 1 & 1 & 1 & 1 & \ldots   \\
   \vdots  & \vdots  & \vdots  & \vdots  & \vdots  & \vdots  & \vdots  & \vdots  & \vdots  & \vdots  & \vdots  & \vdots  & \ddots   \\
\end{matrix} \right),  \quad\left( \begin{matrix}
   1 & 0 & 0 & 0 & 0 & 0 & 0 & 0 & 0 & 0 & 0 & 0 & \ldots   \\
   1 & 1 & 0 & 0 & 0 & 0 & 0 & 0 & 0 & 0 & 0 & 0 & \ldots   \\
   1 & 2 & 1 & 0 & 0 & 0 & 0 & 0 & 0 & 0 & 0 & 0 & \ldots   \\
   1 & 1 & 1 & 1 & 0 & 0 & 0 & 0 & 0 & 0 & 0 & 0 & \ldots   \\
   1 & 0 & 0 & 0 & 1 & 0 & 0 & 0 & 0 & 0 & 0 & 0 & \ldots   \\
   1 & 1 & 0 & 0 & 1 & 1 & 0 & 0 & 0 & 0 & 0 & 0 & \ldots   \\
   1 & 2 & 1 & 0 & 1 & 2 & 1 & 0 & 0 & 0 & 0 & 0 & \ldots   \\
   1 & 1 & 1 & 1 & 1 & 1 & 1 & 1 & 0 & 0 & 0 & 0 & \ldots   \\
   1 & 0 & 0 & 0 & 2 & 0 & 0 & 0 & 1 & 0 & 0 & 0 & \ldots   \\
   1 & 1 & 0 & 0 & 2 & 2 & 0 & 0 & 1 & 1 & 0 & 0 & \ldots   \\
   1 & 2 & 1 & 0 & 2 & 4 & 2 & 0 & 1 & 2 & 1 & 0 & \ldots   \\
   1 & 1 & 1 & 1 & 2 & 2 & 2 & 2 & 1 & 1 & 1 & 1 & \ldots   \\
   \vdots  & \vdots  & \vdots  & \vdots  & \vdots  & \vdots  & \vdots  & \vdots  & \vdots  & \vdots  & \vdots  & \vdots  & \ddots   \\
\end{matrix} \right).$$

If $i<j$, we transform a numerator and denominator of the generalized binomial  coefficient ${{{2^k}n+i\choose {2^k}m+j}_{2}}$ as follows:
$${{b}_{{{2}^{k}}n+i}}!={{b}_{{{2}^{k}}n}}!{{b}_{i}}!={{b}_{{{2}^{k}}n}}{{b}_{{{2}^{k}}n-1}}!{{b}_{i}}!={{b}_{{{2}^{k}}}}{{b}_{n}}{{b}_{{{2}^{k}}\left( n-1 \right)+{{2}^{k}}-1}}!{{b}_{i}}!=$$
$$={{b}_{{{2}^{k}}}}{{b}_{n}}{{b}_{{{2}^{k}}\left( n-1 \right)}}!{{b}_{{{2}^{k}}-1}}!{{b}_{i}}!={{b}_{n}}{{b}_{{{2}^{k}}\left( n-1 \right)}}!{{b}_{{{2}^{k}}+i}}!,$$
where ${{b}_{{{2}^{k}}n}}={{2}^{k}}{{b}_{n}}={{b}_{{{2}^{k}}}}{{b}_{n}}$. For denominator:
$${{b}_{{{2}^{k}}m+j}}!{{b}_{{{2}^{k}}\left( n-m \right)+i-j}}!={{b}_{{{2}^{k}}m}}!{{b}_{j}}!{{b}_{{{2}^{k}}\left( n-m-1 \right)+{{2}^{k}}+i-j}}!={{b}_{{{2}^{k}}m}}!{{b}_{j}}!{{b}_{{{2}^{k}}\left( n-m-1 \right)}}!{{b}_{{{2}^{k}}+i-j}}!.$$
Thus,
$${{\left( \begin{matrix}
   {{2}^{k}}n+i  \\
   {{2}^{k}}m+j  \\
\end{matrix} \right)}_{2}}={{b}_{n}}{{\left( \begin{matrix}
   n-1  \\
   m  \\
\end{matrix} \right)}_{2}}{{\left( \begin{matrix}
   {{2}^{k}}+i  \\
   j  \\
\end{matrix} \right)}_{2}}={{b}_{m+1}}{{\left( \begin{matrix}
   n  \\
   m+1  \\
\end{matrix} \right)}_{2}}{{\left( \begin{matrix}
   {{2}^{k}}+i  \\
   j  \\
\end{matrix} \right)}_{2}},$$  
$$0\le i<{{2}^{k}}, \qquad i<j.$$
For example,
  $$_{\left[ 2 \right]}P-\left( _{\left[ 2 \right]}P\times {}_{0,4}P \right)=\left( \begin{matrix}
   0 & 0 & 0 & 0 & 0 & 0 & 0 & 0 & 0 & 0 & 0 & 0 & \ldots   \\
   0 & 0 & 0 & 0 & 0 & 0 & 0 & 0 & 0 & 0 & 0 & 0 & \ldots   \\
   0 & 0 & 0 & 0 & 0 & 0 & 0 & 0 & 0 & 0 & 0 & 0 & \ldots   \\
   0 & 0 & 0 & 0 & 0 & 0 & 0 & 0 & 0 & 0 & 0 & 0 & \ldots   \\
   0 & 4 & 2 & 4 & 0 & 0 & 0 & 0 & 0 & 0 & 0 & 0 & \ldots   \\
   0 & 0 & 2 & 2 & 0 & 0 & 0 & 0 & 0 & 0 & 0 & 0 & \ldots   \\
   0 & 0 & 0 & 4 & 0 & 0 & 0 & 0 & 0 & 0 & 0 & 0 & \ldots   \\
   0 & 0 & 0 & 0 & 0 & 0 & 0 & 0 & 0 & 0 & 0 & 0 & \ldots   \\
   0 & 8 & 4 & 8 & 0 & 8 & 4 & 8 & 0 & 0 & 0 & 0 & \ldots   \\
   0 & 0 & 4 & 4 & 0 & 0 & 4 & 4 & 0 & 0 & 0 & 0 & \ldots   \\
   0 & 0 & 0 & 8 & 0 & 0 & 0 & 8 & 0 & 0 & 0 & 0 & \ldots   \\
   0 & 0 & 0 & 0 & 0 & 0 & 0 & 0 & 0 & 0 & 0 & 0 & \ldots   \\
   \vdots  & \vdots  & \vdots  & \vdots  & \vdots  & \vdots  & \vdots  & \vdots  & \vdots  & \vdots  & \vdots  & \vdots  & \ddots   \\
\end{matrix} \right).$$

Identities 
$${{\left( \begin{matrix}
   2n+1  \\
   2m+1  \\
\end{matrix} \right)}_{2}}={{\left( \begin{matrix}
   2n+1  \\
   2m  \\
\end{matrix} \right)}_{2}}={{\left( \begin{matrix}
   2n  \\
   2m  \\
\end{matrix} \right)}_{2}}={{\left( \begin{matrix}
   n  \\
   m  \\
\end{matrix} \right)}_{2}},$$

$${{\left( \begin{matrix}
   2n  \\
   2m+1  \\
\end{matrix} \right)}_{2}}=2{{b}_{n}}{{\left( \begin{matrix}
   n-1  \\
   m  \\
\end{matrix} \right)}_{2}}=2{{b}_{m+1}}{{\left( \begin{matrix}
   n  \\
   m+1  \\
\end{matrix} \right)}_{2}}$$
means that rows and columns of the matrix $_{\left[ 2 \right]}P$, denote them ${{u}_{n}}\left( x \right)$ and ${{g}_{n}}\left( x \right)$, form the recurrent sequences:
$${{u}_{2n}}\left( x \right)={{u}_{n}}\left( {{x}^{2}} \right)+2{{b}_{n}}x{{u}_{n-1}}\left( {{x}^{2}} \right),   \qquad{{u}_{2n+1}}\left( x \right)=\left( 1+x \right){{u}_{n}}\left( {{x}^{2}} \right);$$
$${{g}_{2n}}\left( x \right)=\left( 1+x \right){{g}_{n}}\left( {{x}^{2}} \right),    \qquad{{g}_{2n+1}}\left( x \right)=x{{g}_{n}}\left( {{x}^{2}} \right)+2{{b}_{n+1}}{{g}_{n+1}}\left( {{x}^{2}} \right).$$

Turning to generalize, consider the matrix
$$_{\left[ 3 \right]}P={}_{3,3}P\times {}_{{{3,3}^{2}}}P\times {}_{{{3,3}^{3}}}P\times ...\times {}_{{{3,3}^{k}}}P\times ...$$
$$_{\left[ 3 \right]}P=\left( \begin{matrix}
   1 & 0 & 0 & 0 & 0 & 0 & 0 & 0 & 0 & 0 & 0 & 0 & 0 & 0 & 0 & 0 & 0 & 0 & \ldots   \\
   1 & 1 & 0 & 0 & 0 & 0 & 0 & 0 & 0 & 0 & 0 & 0 & 0 & 0 & 0 & 0 & 0 & 0 & \ldots   \\
   1 & 1 & 1 & 0 & 0 & 0 & 0 & 0 & 0 & 0 & 0 & 0 & 0 & 0 & 0 & 0 & 0 & 0 & \ldots   \\
   1 & 3 & 3 & 1 & 0 & 0 & 0 & 0 & 0 & 0 & 0 & 0 & 0 & 0 & 0 & 0 & 0 & 0 & \ldots   \\
   1 & 1 & 3 & 1 & 1 & 0 & 0 & 0 & 0 & 0 & 0 & 0 & 0 & 0 & 0 & 0 & 0 & 0 & \ldots   \\
   1 & 1 & 1 & 1 & 1 & 1 & 0 & 0 & 0 & 0 & 0 & 0 & 0 & 0 & 0 & 0 & 0 & 0 & \ldots   \\
   1 & 3 & 3 & 1 & 3 & 3 & 1 & 0 & 0 & 0 & 0 & 0 & 0 & 0 & 0 & 0 & 0 & 0 & \ldots   \\
   1 & 1 & 3 & 1 & 1 & 3 & 1 & 1 & 0 & 0 & 0 & 0 & 0 & 0 & 0 & 0 & 0 & 0 & \ldots   \\
   1 & 1 & 1 & 1 & 1 & 1 & 1 & 1 & 1 & 0 & 0 & 0 & 0 & 0 & 0 & 0 & 0 & 0 & \ldots   \\
   1 & 9 & 9 & 3 & 9 & 9 & 3 & 9 & 9 & 1 & 0 & 0 & 0 & 0 & 0 & 0 & 0 & 0 & \ldots   \\
   1 & 1 & 9 & 3 & 3 & 9 & 3 & 3 & 9 & 1 & 1 & 0 & 0 & 0 & 0 & 0 & 0 & 0 & \ldots   \\
   1 & 1 & 1 & 3 & 3 & 3 & 3 & 3 & 3 & 1 & 1 & 1 & 0 & 0 & 0 & 0 & 0 & 0 & \ldots   \\
   1 & 3 & 3 & 1 & 9 & 9 & 3 & 9 & 9 & 1 & 3 & 3 & 1 & 0 & 0 & 0 & 0 & 0 & \ldots   \\
   1 & 1 & 3 & 1 & 1 & 9 & 3 & 3 & 9 & 1 & 1 & 3 & 1 & 1 & 0 & 0 & 0 & 0 & \ldots   \\
   1 & 1 & 1 & 1 & 1 & 1 & 3 & 3 & 3 & 1 & 1 & 1 & 1 & 1 & 1 & 0 & 0 & 0 & \ldots   \\
   1 & 3 & 3 & 1 & 3 & 3 & 1 & 9 & 9 & 1 & 3 & 3 & 1 & 3 & 3 & 1 & 1 & 0 & \ldots   \\
   1 & 1 & 3 & 1 & 1 & 3 & 1 & 1 & 9 & 1 & 1 & 3 & 1 & 1 & 3 & 1 & 1 & 1 & \ldots   \\
   1 & 1 & 1 & 1 & 1 & 1 & 1 & 1 & 1 & 1 & 1 & 1 & 1 & 1 & 1 & 1 & 1 & 1 & \ldots   \\
   \vdots  & \vdots  & \vdots  & \vdots  & \vdots  & \vdots  & \vdots  & \vdots  & \vdots  & \vdots  & \vdots  & \vdots  & \vdots  & \vdots  & \vdots  & \vdots  & \vdots  & \vdots  & \ddots   \\
\end{matrix} \right).$$ 
We will use the same notation as in the previous example. The first column of the matrix $_{\left[ 3 \right]}P$ is the generating function of the distribution of divisors ${{3}^{k}}$ in the series of natural numbers:
$$b\left( x \right)=x+{{x}^{2}}+3{{x}^{3}}+{{x}^{4}}+{{x}^{5}}+3{{x}^{6}}+{{x}^{7}}+{{x}^{8}}+9{{x}^{9}}+{{x}^{10}}+{{x}^{11}}+3{{x}^{12}}+$$
$$+{{x}^{13}}+{{x}^{14}}+3{{x}^{15}}+{{x}^{16}}+{{x}^{17}}+9{{x}^{18}}+{{x}^{19}}+{{x}^{20}}+3{{x}^{21}}+{{x}^{22}}+{{x}^{23}}+...,$$
$${{b}_{{{3}^{k}}n}}={{3}^{k}}{{b}_{n}},  \qquad{{b}_{3n+1}}={{b}_{3n+2}}=1,  \qquad{{b}_{{{3}^{k}}n+i}}={{b}_{i}}, \qquad 0<i<{{3}^{k}},$$
  $${{b}_{3n+i}}!={{b}_{3n}}!, \qquad 0\le i<3; \qquad{{b}_{3n-i}}!={{b}_{3\left( n-1 \right)}}!,  \qquad 0<i\le 3,$$
$${{b}_{3n}}!={{3}^{n}}{{b}_{n}}!,  \qquad{{b}_{{{3}^{k}}n}}!={{3}^{\left( \frac{{{3}^{k}}-1}{2} \right)n}}{{b}_{n}}!,$$
$$b\left( x \right)=\frac{\left( 1+x \right)x}{1-{{x}^{3}}}+3b\left( {{x}^{3}} \right)=\sum\limits_{n=0}^{\infty }{\frac{\left( 1+{{x}^{{{3}^{n}}}} \right){{3}^{n}}{{x}^{{{3}^{n}}}}}{1-{{x}^{{{3}^{n+1}}}}}},$$
$${{c}_{3n}}={{c}_{3n+1}}={{c}_{3n+2}}=\frac{{{c}_{n}}}{{{3}^{n}}},$$  
$$c\left( x \right)=\left( 1+x+{{x}^{2}} \right)c\left( \frac{{{x}^{3}}}{3} \right)=\prod\limits_{n=0}^{\infty }{\left( 1+\frac{{{x}^{{{3}^{n}}}}}{{{3}^{{\left( {{3}^{n}}-1 \right)}/{2}\;}}}+\frac{{{x}^{2\left( {{3}^{n}} \right)}}}{{{3}^{{{3}^{n}}-1}}} \right)}.$$
If
$$\frac{{{b}_{n}}!}{{{b}_{m}}!{{b}_{n-m}}!}={{\left( \begin{matrix}
   n  \\
   m  \\
\end{matrix} \right)}_{3}},$$
then
$${{\left( \begin{matrix}
   {{3}^{k}}n+i  \\
   {{3}^{k}}m+j  \\
\end{matrix} \right)}_{3}}={{\left( \begin{matrix}
   n  \\
   m  \\
\end{matrix} \right)}_{3}}{{\left( \begin{matrix}
   i  \\
   j  \\
\end{matrix} \right)}_{3}}, \qquad0\le i<{{3}^{k}}, \qquad i\ge j;$$
$${{\left( \begin{matrix}
   {{3}^{k}}n+i  \\
   {{3}^{k}}m+j  \\
\end{matrix} \right)}_{3}}={{b}_{n}}{{\left( \begin{matrix}
   n-1  \\
   m  \\
\end{matrix} \right)}_{3}}{{\left( \begin{matrix}
   {{3}^{k}}+i  \\
   j  \\
\end{matrix} \right)}_{3}}={{b}_{m+1}}{{\left( \begin{matrix}
   n  \\
   m+1  \\
\end{matrix} \right)}_{3}}{{\left( \begin{matrix}
   {{3}^{k}}+i  \\
   j  \\
\end{matrix} \right)}_{3}},$$ 
$$0\le i<{{3}^{k}}, \qquad i<j.$$
From
$${{\left( \begin{matrix}
   3n+i  \\
   3m+j  \\
\end{matrix} \right)}_{3}}={{\left( \begin{matrix}
   n  \\
   m  \\
\end{matrix} \right)}_{3}}, \qquad i\ge j; \qquad=3{{b}_{n}}{{\left( \begin{matrix}
   n-1  \\
   m  \\
\end{matrix} \right)}_{3}}=3{{b}_{m+1}}{{\left( \begin{matrix}
   n  \\
   m+1  \\
\end{matrix} \right)}_{3}}, \qquad i<j,$$
we see that if
$$\left[ n,\to  \right]{}_{\left[ 3 \right]}P={{u}_{n}}\left( x \right), \qquad\left[ \uparrow ,n \right]{}_{\left[ 3 \right]}P={{g}_{n}}\left( x \right),$$
then
$${{u}_{3n}}\left( x \right)={{u}_{n}}\left( {{x}^{3}} \right)+3{{b}_{n}}\left( 1+x \right)x{{u}_{n-1}}\left( {{x}^{3}} \right),$$ 
$${{u}_{3n+1}}\left( x \right)=\left( 1+x \right){{u}_{n}}\left( {{x}^{3}} \right)+3{{b}_{n}}{{x}^{2}}{{u}_{n-1}}\left( {{x}^{3}} \right),$$  
$${{u}_{3n+2}}\left( x \right)=\left( 1+x+{{x}^{2}} \right){{u}_{n}}\left( {{x}^{3}} \right);$$
$${{g}_{3n}}\left( x \right)=\left( 1+x+{{x}^{2}} \right){{g}_{n}}\left( {{x}^{3}} \right),$$ 
$${{g}_{3n+1}}\left( x \right)=\left( 1+x \right)x{{g}_{n}}\left( {{x}^{3}} \right)+3{{b}_{n+1}}{{g}_{n+1}}\left( {{x}^{3}} \right),$$
$${{g}_{3n+2}}\left( x \right)={{x}^{2}}{{g}_{n}}\left( {{x}^{3}} \right)+3{{b}_{n+1}}\left( 1+x \right){{g}_{n+1}}\left( {{x}^{3}} \right).$$

Introduce the notation
$${{w}_{m}}\left( x \right)=\sum\limits_{n=0}^{m}{{{x}^{n}}},  \qquad{{w}_{-1}}\left( x \right)=0.$$
For the matrix
$$_{\left[ q \right]}P={}_{q,q}P\times {}_{q,{{q}^{2}}}P\times {}_{q,{{q}^{3}}}P\times ...\times {}_{q,{{q}^{k}}}P\times ...,$$
$q=2$, $3$, $4$, … , we have (in the same notation as in the previous examples):
$${{b}_{{{q}^{k}}n}}={{q}^{k}}{{b}_{n}},  \qquad{{b}_{qn+i}}=1, \qquad0<i<q, \qquad{{b}_{{{q}^{k}}n+i}}={{b}_{i}}, \qquad0<i<{{q}^{k}}.$$
$${{b}_{{{q}^{k}}n}}!={{q}^{\left( \frac{{{q}^{k}}-1}{q-1} \right)n}}{{b}_{n}}!,$$
$$b\left( x \right)=\frac{{{w}_{q-2}}\left( x \right)x}{1-{{x}^{q}}}+qb\left( {{x}^{q}} \right)=\sum\limits_{n=0}^{\infty }{\frac{{{w}_{q-2}}\left( {{x}^{{{q}^{n}}}} \right){{q}^{n}}{{x}^{{{q}^{n}}}}}{1-{{x}^{{{q}^{n+1}}}}}};$$
$${{c}_{qn+i}}=\frac{{{c}_{n}}}{{{q}^{n}}},  \qquad0\le i<q,$$
$$c\left( x \right)={{w}_{q-1}}\left( x \right)c\left( \frac{{{x}^{q}}}{q} \right)=\prod\limits_{n=0}^{\infty }{{{w}_{q-1}}}\left( {{{x}^{{{q}^{n}}}}}/{{{q}^{\frac{{{q}^{n}}-1}{q-1}}}}\; \right);$$
$${{\left( \begin{matrix}
   {{q}^{k}}n+i  \\
   {{q}^{k}}m+j  \\
\end{matrix} \right)}_{q}}={{\left( \begin{matrix}
   n  \\
   m  \\
\end{matrix} \right)}_{q}}{{\left( \begin{matrix}
   i  \\
   j  \\
\end{matrix} \right)}_{q}},  \qquad0\le i<{{q}^{k}}, \qquad i\ge j,$$

$${{\left( \begin{matrix}
   {{q}^{k}}n+i  \\
   {{q}^{k}}m+j  \\
\end{matrix} \right)}_{q}}={{b}_{n}}{{\left( \begin{matrix}
   n-1  \\
   m  \\
\end{matrix} \right)}_{q}}{{\left( \begin{matrix}
   {{q}^{k}}+i  \\
   j  \\
\end{matrix} \right)}_{q}}={{b}_{m+1}}{{\left( \begin{matrix}
   n  \\
   m+1  \\
\end{matrix} \right)}_{q}}{{\left( \begin{matrix}
   {{q}^{k}}+i  \\
   j  \\
\end{matrix} \right)}_{q}}, \qquad i<j;$$
$${{u}_{qn+m}}\left( x \right)={{w}_{m}}\left( x \right){{u}_{n}}\left( {{x}^{q}} \right)+q{{b}_{n}}{{x}^{m+1}}{{w}_{q-2-m}}\left( x \right){{u}_{n-1}}\left( {{x}^{q}} \right),$$
$${{g}_{qn+m}}\left( x \right)={{x}^{m}}{{w}_{q-1-m}}\left( x \right){{g}_{n}}\left( {{x}^{q}} \right)+q{{b}_{n+1}}{{w}_{m-1}}\left( x \right){{g}_{n+1}}\left( {{x}^{q}} \right), \quad 0\le m<q.$$

According to the definition of the matrix $_{\left[ q \right]}P$,
$$P={}_{\left[ 2 \right]}P\times {}_{\left[ 3 \right]}P\times {}_{\left[ 5 \right]}P\times ...\times {}_{\left[ p \right]}P\times ...,$$
where $p$ is a member of the sequence of prime numbers. Respectively, the exponential series is the Hadamard product of the series
$$\prod\limits_{n=0}^{\infty }{{{w}_{p-1}}}\left( {{{x}^{{{p}^{n}}}}}/{{{p}^{\frac{{{p}^{n}}-1}{p-1}}}}\; \right).$$

Generalization of the matrix $_{\left[ q \right]}P$ is the matrix
$$_{\left[ \varphi ,q \right]}P={}_{\varphi ,q}P\times {}_{\varphi ,{{q}^{2}}}P\times {}_{\varphi ,{{q}^{3}}}P\times ...\times {}_{\varphi ,{{q}^{k}}}P\times ...,$$ 
${}_{\left[ q,q \right]}P={}_{\left[ q \right]}P$. If $\varphi \ne 0$, identities for the matrix $_{\left[ \varphi ,q \right]}P$ can be obtained from the identities for the matrix $_{\left[ q \right]}P$ by replacing ${{b}_{{{q}^{k}}n}}={{q}^{k}}{{b}_{n}}$, ${{c}_{qn+i}}={{{c}_{n}}}/{{{q}^{n}}}\;$ on ${{b}_{{{q}^{k}}n}}={{\varphi }^{k}}{{b}_{n}}$, ${{c}_{qn+i}}={{{c}_{n}}}/{{{\varphi }^{n}}}\;$. If $q$ is fixed, matrices $_{\left[ \varphi ,q \right]}P$ form the group:
$${}_{\left[ \varphi ,q \right]}P\times {}_{\left[ \beta ,q \right]}P={}_{\left[ \varphi \beta ,q \right]}P,$$
where $_{\left[ 1,q \right]}P={}_{1,q}P={{P}_{{{\left( 1-x \right)}^{-1}}}}$ is the identity element of  the group of generalized Pascal matrices.
If $\varphi =0$, we get zero generalized Pascal matrix
\[_{\left[ 0,q \right]}P={}_{0,q}P\times {}_{0,{{q}^{2}}}P\times {}_{0,{{q}^{3}}}P\times ...\times {}_{0,{{q}^{k}}}P\times ...\]
For example (Pascal triangle modulo 2),
$$_{\left[ 0,2 \right]}P=\left( \begin{matrix}
   1 & 0 & 0 & 0 & 0 & 0 & 0 & 0 & 0 & 0 & 0 & 0 & 0 & 0 & 0 & 0 & \ldots   \\
   1 & 1 & 0 & 0 & 0 & 0 & 0 & 0 & 0 & 0 & 0 & 0 & 0 & 0 & 0 & 0 & \ldots   \\
   1 & 0 & 1 & 0 & 0 & 0 & 0 & 0 & 0 & 0 & 0 & 0 & 0 & 0 & 0 & 0 & \ldots   \\
   1 & 1 & 1 & 1 & 0 & 0 & 0 & 0 & 0 & 0 & 0 & 0 & 0 & 0 & 0 & 0 & \ldots   \\
   1 & 0 & 0 & 0 & 1 & 0 & 0 & 0 & 0 & 0 & 0 & 0 & 0 & 0 & 0 & 0 & \ldots   \\
   1 & 1 & 0 & 0 & 1 & 1 & 0 & 0 & 0 & 0 & 0 & 0 & 0 & 0 & 0 & 0 & \ldots   \\
   1 & 0 & 1 & 0 & 1 & 0 & 1 & 0 & 0 & 0 & 0 & 0 & 0 & 0 & 0 & 0 & \ldots   \\
   1 & 1 & 1 & 1 & 1 & 1 & 1 & 1 & 0 & 0 & 0 & 0 & 0 & 0 & 0 & 0 & \ldots   \\
   1 & 0 & 0 & 0 & 0 & 0 & 0 & 0 & 1 & 0 & 0 & 0 & 0 & 0 & 0 & 0 & \ldots   \\
   1 & 1 & 0 & 0 & 0 & 0 & 0 & 0 & 1 & 1 & 0 & 0 & 0 & 0 & 0 & 0 & \ldots   \\
   1 & 0 & 1 & 0 & 0 & 0 & 0 & 0 & 1 & 0 & 1 & 0 & 0 & 0 & 0 & 0 & \ldots   \\
   1 & 1 & 1 & 1 & 0 & 0 & 0 & 0 & 1 & 1 & 1 & 1 & 0 & 0 & 0 & 0 & \ldots   \\
   1 & 0 & 0 & 0 & 1 & 0 & 0 & 0 & 1 & 0 & 0 & 0 & 1 & 0 & 0 & 0 & \ldots   \\
   1 & 1 & 0 & 0 & 1 & 1 & 0 & 0 & 1 & 1 & 0 & 0 & 1 & 1 & 0 & 0 & \ldots   \\
   1 & 0 & 1 & 0 & 1 & 0 & 1 & 0 & 1 & 0 & 1 & 0 & 1 & 0 & 1 & 0 & \ldots   \\
   1 & 1 & 1 & 1 & 1 & 1 & 1 & 1 & 1 & 1 & 1 & 1 & 1 & 1 & 1 & 1 & \ldots   \\
   \vdots  & \vdots  & \vdots  & \vdots  & \vdots  & \vdots  & \vdots  & \vdots  & \vdots  & \vdots  & \vdots  & \vdots  & \vdots  & \vdots  & \vdots  & \vdots  & \ddots   \\
\end{matrix} \right).$$
$${{\left( {}_{\left[ 0,q \right]}P \right)}_{n,m}}=1, \quad n\left( \bmod {{q}^{k}} \right)\ge m\left( \bmod {{q}^{k}} \right);
\quad=0,\quad n\left( \bmod {{q}^{k}} \right)<m\left( \bmod {{q}^{k}} \right),\quad k=1, 2, … $$

\section{Fractal zero generalized Pascal matrices}

Denote
$${{\left( _{\left[ 0,q \right]}P \right)}_{n,m}}={{\left( \begin{matrix}
   n  \\
   m  \\
\end{matrix} \right)}_{0,q}}.$$
{\bfseries Theorem.}
$${{\left( \begin{matrix}
   {{q}^{k}}n+i  \\
   {{q}^{k}}m+j  \\
\end{matrix} \right)}_{0,q}}={{\left( \begin{matrix}
   n  \\
   m  \\
\end{matrix} \right)}_{0,q}}{{\left( \begin{matrix}
   i  \\
   j  \\
\end{matrix} \right)}_{0,q}},  \qquad 0\le i,j<{{q}^{k}}.\eqno(3)$$ 
{\bfseries Proof.} By definition, if $n\left( \bmod {{q}^{k}} \right)<m\left( \bmod {{q}^{k}} \right)$ for some value of $k$, then ${{n\choose m}_{0,q}}=0$. We represent the numbers $n$, $m$ in the form 
$$n=\sum\limits_{i=0}^{\infty }{{{n}_{i}}}{{q}^{i}},  \qquad m=\sum\limits_{i=0}^{\infty }{{{m}_{i}}{{q}^{i}}}, \qquad 0\le {{n}_{i}},{{m}_{i}}<q.$$
Then
$$n\left( \bmod {{q}^{k}} \right)=\sum\limits_{i=0}^{k-1}{{{n}_{i}}}{{q}^{i}},  \qquad m\left( \bmod {{q}^{k}} \right)=\sum\limits_{i=0}^{k-1}{{{m}_{i}}{{q}^{i}}}.$$
If $n\left( \bmod {{q}^{k}} \right)<m\left( \bmod {{q}^{k}} \right)$, then ${{n}_{i}}<{{m}_{i}}$ at least for one $i$. Since ${{n\choose m}_{0,q}}=1$, if ${{m}_{i}}\ge {{n}_{i}}$, then true the identity
$${{\left( \begin{matrix}
   n  \\
   m  \\
\end{matrix} \right)}_{0,q}}=\prod\limits_{i=0}^{\infty }{{{\left( \begin{matrix}
   {{n}_{i}}  \\
   {{m}_{i}}  \\
\end{matrix} \right)}_{0,q}}}.\eqno(4)$$
It remains to note that if
$$n=\sum\limits_{i=0}^{\infty }{{{n}_{i}}{{q}^{i}}=}{{q}^{k}}s+j, \qquad 0\le j<{{q}^{k}},$$ 
then
$$j=\sum\limits_{i=0}^{k-1}{{{n}_{i}}}{{q}^{i}}, \qquad s=\sum\limits_{i=0}^{\infty }{{{n}_{i+k}}{{q}^{i}}}.$$

 In the works [7], [8], [9] matrices $_{\left[ 0,q \right]}P$ are called generalized Sierpinski matrices. The property (3) is represented as 
$${{S}_{q}}={{S}_{q,k}}\otimes {{S}_{q,k}}\otimes {{S}_{q,k}}\otimes ... ,   \qquad k=1, 2 , … ,$$
where ${{S}_{q}}{{=}_{\left[ 0,q \right]}}P$, ${{S}_{q,k}}$ is the matrix consisting of  ${{q}^{k}}$ first rows of the matrix ${{S}_{q}}$, $\otimes $ denotes the Kronecker multiplication. Matrices ${{S}_{q}}$ have a generalizations, one of which can be represented as
$${{S}_{q}}\left( a\left( x \right) \right)={{S}_{q,k}}\left( a\left( x \right) \right)\otimes {{S}_{q,k}}\left( a\left( x \right) \right)\otimes {{S}_{q,k}}\left( a\left( x \right) \right)\otimes ... ,$$
where ${{S}_{q,1}}\left( a\left( x \right) \right)$ is the matrix consisting of $q$ first rows of the matrix $A$:
$$\left[ \uparrow ,n \right]A={{x}^{n}}a\left( x \right), \qquad{{\left( A \right)}_{n,m}}={{a}_{n-m}},$$
${{S}_{q,k}}\left( a\left( x \right) \right)$  is the matrix consisting of ${{q}^{k}}$ first rows of the matrix ${{S}_{q}}\left( a\left( x \right) \right)$. Then
$${{S}_{q}}\left( a\left( x \right) \right){{S}_{q}}\left( b\left( x \right) \right)={{S}_{q}}\left( a\left( x \right)b\left( x \right) \right).$$
The results of these works have been developed in [10], [11]. Introduced matrices
$${{T}^{\left( q \right)}}=T_{k}^{\left( q \right)}\otimes T_{k}^{\left( q \right)}\otimes T_{k}^{\left( q \right)}\otimes ... ,  \qquad k=1, 2 , … ,$$
where $T_{1}^{\left( q \right)}$ is the matrix consisting of $q$ first rows of the Pascal matrix, $T_{k}^{\left( q \right)}$ is the matrix consisting of ${{q}^{k}}$ first rows of the matrix ${{T}^{\left( q \right)}}$. For example,
$${{T}^{\left( 3 \right)}}=\left( \begin{matrix}
   1 & 0 & 0 & 0 & 0 & 0 & 0 & 0 & 0 & \ldots   \\
   1 & 1 & 0 & 0 & 0 & 0 & 0 & 0 & 0 & \ldots   \\
   1 & 2 & 1 & 0 & 0 & 0 & 0 & 0 & 0 & \ldots   \\
   1 & 0 & 0 & 1 & 0 & 0 & 0 & 0 & 0 & \ldots   \\
   1 & 1 & 0 & 1 & 1 & 0 & 0 & 0 & 0 & \ldots   \\
   1 & 2 & 1 & 1 & 2 & 1 & 0 & 0 & 0 & \ldots   \\
   1 & 0 & 0 & 2 & 0 & 0 & 1 & 0 & 0 & \ldots   \\
   1 & 1 & 0 & 2 & 2 & 0 & 1 & 1 & 0 & \ldots   \\
   1 & 2 & 1 & 2 & 4 & 2 & 1 & 2 & 1 & \ldots   \\
   \vdots  & \vdots  & \vdots  & \vdots  & \vdots  & \vdots  & \vdots  & \vdots  & \vdots  & \ddots   \\
\end{matrix} \right).$$
Denote
$${{\left( {{T}^{\left( q \right)}} \right)}_{n,m}}={{\left( \begin{matrix}
   n  \\
   m  \\
\end{matrix} \right)}_{q}}.$$
Then
$${{\left( \begin{matrix}
   n  \\
   m  \\
\end{matrix} \right)}_{q}}=\prod\limits_{i=0}^{\infty }{\left( \begin{matrix}
   {{n}_{i}}  \\
   {{m}_{i}}  \\
\end{matrix} \right)},  \quad n=\sum\limits_{i=0}^{\infty }{{{n}_{i}}}{{q}^{i}},  \quad m=\sum\limits_{i=0}^{\infty }{{{m}_{i}}}{{q}^{i}},  \quad 0\le {{n}_{i}},{{m}_{i}}<q,$$
$$\sum\limits_{m=0}^{n}{{{\left( \begin{matrix}
   n  \\
   m  \\
\end{matrix} \right)}_{q}}}{{x}^{m}}=\prod\limits_{i=0}^{\infty }{{{\left( 1+{{x}^{{{q}^{i}}}} \right)}^{{{n}_{i}}}}}.$$

These results become more transparent if we use algebra of the matrices $\left( a\left( x \right)|q \right)$:
$${{\left( \left( a\left( x \right)|q \right) \right)}_{n,m}}={{a}_{n-m}}{{\left( \begin{matrix}
   n  \\
   m  \\
\end{matrix} \right)}_{0,q}},$$
$$\left( a\left( x \right)|q \right)\left( b\left( x \right)|q \right)=\left( a\left( x \right)\circ b\left( x \right)|q \right),$$
$$\left[ {{x}^{n}} \right]a\left( x \right)\circ b\left( x \right)=\sum\limits_{m=0}^{n}{{{\left( \begin{matrix}
   n  \\
   m  \\
\end{matrix} \right)}_{0,q}}}{{a}_{m}}{{b}_{n-m}}.$$
Denote
$$\left( \left( \sum\limits_{n=0}^{{{q}^{k}}-1}{{{b}_{n}}{{x}^{n}}} \right)a\left( {{x}^{{{q}^{k}}}} \right)|q \right)=\left( a\left( x \right),b\left( x \right)|q,k \right).$$
For example,
$$\left( a\left( x \right),b\left( x \right)|2,1 \right)=\left( \begin{matrix}
   {{a}_{0}}{{b}_{0}} & 0 & 0 & 0 & 0 & 0 & 0 & 0 & \ldots   \\
   {{a}_{0}}{{b}_{1}} & {{a}_{0}}{{b}_{0}} & 0 & 0 & 0 & 0 & 0 & 0 & \ldots   \\
   {{a}_{1}}{{b}_{0}} & 0 & {{a}_{0}}{{b}_{0}} & 0 & 0 & 0 & 0 & 0 & \ldots   \\
   {{a}_{1}}{{b}_{1}} & {{a}_{1}}{{b}_{0}} & {{a}_{0}}{{b}_{1}} & {{a}_{0}}{{b}_{0}} & 0 & 0 & 0 & 0 & \ldots   \\
   {{a}_{2}}{{b}_{0}} & 0 & 0 & 0 & {{a}_{0}}{{b}_{0}} & 0 & 0 & 0 & \ldots   \\
   {{a}_{2}}{{b}_{1}} & {{a}_{2}}{{b}_{0}} & 0 & 0 & {{a}_{0}}{{b}_{1}} & {{a}_{0}}{{b}_{0}} & 0 & 0 & \ldots   \\
   {{a}_{3}}{{b}_{0}} & 0 & {{a}_{2}}{{b}_{0}} & 0 & {{a}_{1}}{{b}_{0}} & 0 & {{a}_{0}}{{b}_{0}} & 0 & \ldots   \\
   {{a}_{3}}{{b}_{1}} & {{a}_{3}}{{b}_{0}} & {{a}_{2}}{{b}_{1}} & {{a}_{2}}{{b}_{0}} & {{a}_{1}}{{b}_{1}} & {{a}_{1}}{{b}_{0}} & {{a}_{0}}{{b}_{1}} & {{a}_{0}}{{b}_{0}} & \ldots   \\
   \vdots  & \vdots  & \vdots  & \vdots  & \vdots  & \vdots  & \vdots  & \vdots  & \ddots   \\
\end{matrix} \right),$$
$$\left( a\left( x \right),b\left( x \right)|2,2 \right)=\left( \begin{matrix}
   {{a}_{0}}{{b}_{0}} & 0 & 0 & 0 & 0 & 0 & 0 & 0 & \ldots   \\
   {{a}_{0}}{{b}_{1}} & {{a}_{0}}{{b}_{0}} & 0 & 0 & 0 & 0 & 0 & 0 & \ldots   \\
   {{a}_{0}}{{b}_{2}} & 0 & {{a}_{0}}{{b}_{0}} & 0 & 0 & 0 & 0 & 0 & \ldots   \\
   {{a}_{0}}{{b}_{3}} & {{a}_{0}}{{b}_{2}} & {{a}_{0}}{{b}_{1}} & {{a}_{0}}{{b}_{0}} & 0 & 0 & 0 & 0 & \ldots   \\
   {{a}_{1}}{{b}_{0}} & 0 & 0 & 0 & {{a}_{0}}{{b}_{0}} & 0 & 0 & 0 & \ldots   \\
   {{a}_{1}}{{b}_{1}} & {{a}_{1}}{{b}_{0}} & 0 & 0 & {{a}_{0}}{{b}_{1}} & {{a}_{0}}{{b}_{0}} & 0 & 0 & \ldots   \\
   {{a}_{1}}{{b}_{2}} & 0 & {{a}_{1}}{{b}_{0}} & 0 & {{a}_{0}}{{b}_{2}} & 0 & {{a}_{0}}{{b}_{0}} & 0 & \ldots   \\
   {{a}_{1}}{{b}_{3}} & {{a}_{1}}{{b}_{2}} & {{a}_{1}}{{b}_{1}} & {{a}_{1}}{{b}_{0}} & {{a}_{0}}{{b}_{3}} & {{a}_{0}}{{b}_{2}} & {{a}_{0}}{{b}_{1}} & {{a}_{0}}{{b}_{0}} & \ldots   \\
   \vdots  & \vdots  & \vdots  & \vdots  & \vdots  & \vdots  & \vdots  & \vdots  & \ddots   \\
\end{matrix} \right).$$
Since
$$\left[ {{x}^{{{q}^{k}}n+i}} \right]\left( \sum\limits_{n=0}^{{{q}^{k}}-1}{{{b}_{n}}{{x}^{n}}} \right)a\left( {{x}^{{{q}^{k}}}} \right)={{a}_{n}}{{b}_{i}}, \quad 0\le i<{{q}^{k}},$$
then $\left( {{q}^{k}}n+i,{{q}^{k}}m+j \right)$-th element of the matrix $\left( a\left( x \right),b\left( x \right)|q,k \right)$ is equal to
$${{a}_{n-m}}{{b}_{i-j}}{{\left( \begin{matrix}
   {{q}^{k}}n+i  \\
   {{q}^{k}}m+j  \\
\end{matrix} \right)}_{0,q}}={{a}_{n-m}}{{\left( \begin{matrix}
   n  \\
   m  \\
\end{matrix} \right)}_{0,q}}{{b}_{i-j}}{{\left( \begin{matrix}
   i  \\
   j  \\
\end{matrix} \right)}_{0,q}},\qquad 0\le i,j<{{q}^{k}}.$$
Thus, $\left( a\left( x \right),b\left( x \right)|q,k \right)$ is the block matrix, $\left( n,m \right)$-th block of which is the matrix 
$${{a}_{n-m}}{{\left( \begin{matrix}
   n  \\
   m  \\
\end{matrix} \right)}_{0,q}}{{\left( b\left( x \right)|q \right)}_{{{q}^{k}}}},$$
where ${{\left( b\left( x \right)|q \right)}_{{{q}^{k}}}}$ is the matrix consisting of ${{q}^{k}}$ first rows of the matrix $\left( b\left( x \right)|q \right)$. Hence,
$$\left( a\left( x \right),b\left( x \right)|q,k \right)\left( c\left( x \right),d\left( x \right)|q,k \right)=\left( a\left( x \right)\circ c\left( x \right),b\left( x \right)\circ d\left( x \right)|q,k \right).$$

If
$$\left( a\left( x \right)|q \right)=\left( a\left( x \right),a\left( x \right)|q,1 \right),$$
as in the case of the matrix $_{\left[ 0,q \right]}P$, then $\left( a\left( x \right)|q \right)=\left( a\left( x \right),a\left( x \right)|q,k \right)$ for all $k=1, 2, … $ :
$$a\left( x \right)=\left( \sum\limits_{n=0}^{q-1}{{{a}_{n}}{{x}^{n}}} \right)a\left( {{x}^{q}} \right)=\left( \sum\limits_{n=0}^{{{q}^{k}}-1}{{{a}_{n}}{{x}^{n}}} \right)a\left( {{x}^{{{q}^{k}}}} \right)=\prod\limits_{m=0}^{\infty }{\left( \sum\limits_{n=0}^{{{q}^{k}}-1}{{{a}_{n}}{{x}^{n{{q}^{mk}}}}} \right)},$$
$${{a}_{0}}=1,  \qquad{{a}_{{{q}^{k}}n+i}}={{a}_{n}}{{a}_{i}},  \qquad0\le i<{{q}^{k}},$$
$${{\left( \left( a\left( x \right)|q \right) \right)}_{{{q}^{k}}n+i,{{q}^{k}}m+j}}={{\left( \left( a\left( x \right)|q \right) \right)}_{n,m}}{{\left( \left( a\left( x \right)|q \right) \right)}_{i,j}}.\eqno(5)$$
Hence
$${{a}_{n}}=\prod\limits_{i=0}^{\infty }{{{a}_{{{n}_{i}}}}}, \quad n=\sum\limits_{i=0}^{\infty }{{{n}_{i}}}{{q}^{i}}={{n}_{0}}+q\left( {{n}_{1}}+q\left( {{n}_{2}}+... \right) \right),  \quad0\le {{n}_{i}}<q.$$
In view of the identity (4),
$${{a}_{n-m}}{{\left( \begin{matrix}
   n  \\
   m  \\
\end{matrix} \right)}_{0,q}}=\prod\limits_{i=0}^{\infty }{{{a}_{{{n}_{i}}-{{m}_{i}}}}},  \quad n=\sum\limits_{i=0}^{\infty }{{{n}_{i}}}{{q}^{i}},  \quad m=\sum\limits_{i=0}^{\infty }{{{m}_{i}}}{{q}^{i}},$$
 if ${{n}_{i}}\ge {{m}_{i}}$ for all $i$, and ${{a}_{n-m}}{{n\choose m}_{0,q}}=0$ in other case. I.e.
$${{\left( \left( a\left( x \right)|q \right) \right)}_{n,m}}=\prod\limits_{i=0}^{\infty }{{{\left( \left( a\left( x \right)|q \right) \right)}_{{{n}_{i}},{{m}_{i}}}}}.\eqno(6)$$   
From (5), (6) we see that rows of the matrix $\left( a\left( x \right)|q \right)$, denote them ${{u}_{n}}\left( x \right)$, form the recurrent sequences:
$${{u}_{n}}\left( x \right)=\sum\limits_{m=0}^{n}{{{a}_{n-m}}{{x}^{m}}}, \quad 0\le n<q;  \quad{{u}_{qn+i}}\left( x \right)={{u}_{n}}\left( {{x}^{q}} \right){{u}_{i}}\left( x \right), \quad0\le i<q;$$
or
$${{u}_{n}}\left( x \right)=\prod\limits_{i=0}^{\infty }{{{u}_{{{n}_{i}}}}}\left( {{x}^{{{q}^{i}}}} \right),  \quad n=\sum\limits_{i=0}^{\infty }{{{n}_{i}}}{{q}^{i}},  \quad0\le {{n}_{i}}<q.$$
If
$$\left( a\left( x \right)|q \right)\left( b\left( x \right)|q \right)=\left( a\left( x \right)\circ b\left( x \right)|q \right),$$
where
$$a\left( x \right)=\left( \sum\limits_{n=0}^{q-1}{{{a}_{n}}{{x}^{n}}} \right)a\left( {{x}^{q}} \right), \qquad b\left( x \right)=\left( \sum\limits_{n=0}^{q-1}{{{b}_{n}}{{x}^{n}}} \right)b\left( {{x}^{q}} \right),$$
then
$$\left[ {{x}^{n}} \right]a\left( x \right)\circ b\left( x \right)=\prod\limits_{i=0}^{\infty }{{{c}_{{{n}_{i}}}}},  \quad{{c}_{{{n}_{i}}}}=\left[ {{x}^{{{n}_{i}}}} \right]a\left( x \right)b\left( x \right),  \quad n=\sum\limits_{i=0}^{\infty }{{{n}_{i}}}{{q}^{i}}, \quad0\le {{n}_{i}}<q.$$
For example,
$$_{\left[ 0,2 \right]}{{P}^{2}}={{\left( \frac{1}{1-x}|2 \right)}^{2}}=\left( \begin{matrix}
   1 & 0 & 0 & 0 & 0 & 0 & 0 & 0 & 0 & 0 & 0 & 0 & 0 & 0 & 0 & 0 & \ldots   \\
   2 & 1 & 0 & 0 & 0 & 0 & 0 & 0 & 0 & 0 & 0 & 0 & 0 & 0 & 0 & 0 & \ldots   \\
   2 & 0 & 1 & 0 & 0 & 0 & 0 & 0 & 0 & 0 & 0 & 0 & 0 & 0 & 0 & 0 & \ldots   \\
   4 & 2 & 2 & 1 & 0 & 0 & 0 & 0 & 0 & 0 & 0 & 0 & 0 & 0 & 0 & 0 & \ldots   \\
   2 & 0 & 0 & 0 & 1 & 0 & 0 & 0 & 0 & 0 & 0 & 0 & 0 & 0 & 0 & 0 & \ldots   \\
   4 & 2 & 0 & 0 & 2 & 1 & 0 & 0 & 0 & 0 & 0 & 0 & 0 & 0 & 0 & 0 & \ldots   \\
   4 & 0 & 2 & 0 & 2 & 0 & 1 & 0 & 0 & 0 & 0 & 0 & 0 & 0 & 0 & 0 & \ldots   \\
   8 & 4 & 4 & 2 & 4 & 2 & 2 & 1 & 0 & 0 & 0 & 0 & 0 & 0 & 0 & 0 & \ldots   \\
   2 & 0 & 0 & 0 & 0 & 0 & 0 & 0 & 1 & 0 & 0 & 0 & 0 & 0 & 0 & 0 & \ldots   \\
   4 & 2 & 0 & 0 & 0 & 0 & 0 & 0 & 2 & 1 & 0 & 0 & 0 & 0 & 0 & 0 & \ldots   \\
   4 & 0 & 2 & 0 & 0 & 0 & 0 & 0 & 2 & 0 & 1 & 0 & 0 & 0 & 0 & 0 & \ldots   \\
   8 & 4 & 4 & 2 & 0 & 0 & 0 & 0 & 4 & 2 & 2 & 1 & 0 & 0 & 0 & 0 & \ldots   \\
   4 & 0 & 0 & 0 & 2 & 0 & 0 & 0 & 2 & 0 & 0 & 0 & 1 & 0 & 0 & 0 & \ldots   \\
   8 & 4 & 0 & 0 & 4 & 2 & 0 & 0 & 4 & 2 & 0 & 0 & 2 & 1 & 0 & 0 & \ldots   \\
   8 & 0 & 4 & 0 & 4 & 0 & 2 & 0 & 4 & 0 & 2 & 0 & 2 & 0 & 1 & 0 & \ldots   \\
   16 & 8 & 8 & 4 & 8 & 4 & 4 & 2 & 8 & 4 & 4 & 2 & 4 & 2 & 2 & 1 & \ldots   \\
   \vdots  & \vdots  & \vdots  & \vdots  & \vdots  & \vdots  & \vdots  & \vdots  & \vdots  & \vdots  & \vdots  & \vdots  & \vdots  & \vdots  & \vdots  & \vdots  & \ddots   \\
\end{matrix} \right),$$
$$\left[ {{x}^{n}} \right]{{\left( 1-x \right)}^{-1}}\circ {{\left( 1-x \right)}^{-1}}=\prod\limits_{i=0}^{\infty }{{{2}^{{{n}_{i}}}}},  \quad n=\sum\limits_{i=0}^{\infty }{{{n}_{i}}}{{2}^{i}}, \quad 0\le {{n}_{i}}<2,$$ 
$$\left[ n,\to  \right]{}_{\left[ 0,2 \right]}{{P}^{2}}=\prod\limits_{i=0}^{\infty }{{{\left( 2+{{x}^{{{2}^{i}}}} \right)}^{{{n}_{i}}}}}.$$

Let
$$c\left( x \right)=\left( \sum\limits_{n=0}^{q-1}{{{c}_{n}}{{x}^{n}}} \right)c\left( {{x}^{q}} \right).$$
Then
$${{\left( {{P}_{c\left( x \right)}} \right)}_{{{q}^{k}}n+i,{{q}^{k}}m+j}}=\frac{{{c}_{{{q}^{k}}n+i}}{{c}_{{{q}^{k}}\left( n-m \right)+i-j}}}{{{c}_{{{q}^{k}}m+j}}}=\frac{{{c}_{n}}{{c}_{i}}{{c}_{n-m}}{{c}_{i-j}}}{{{c}_{m}}{{c}_{j}}}={{\left( {{P}_{c\left( x \right)}} \right)}_{n,m}}{{\left( {{P}_{c\left( x \right)}} \right)}_{i,j}},$$ 
$0\le i,j<{{q}^{k}}$, $i\ge j$. Hence, matrices ${{P}_{c\left( x \right)}}\times {}_{\left[ 0,q \right]}P$ and $\left( c\left( x \right)|q \right)$ have the same fractal properties. Applied to the matrix ${{T}^{\left( q \right)}}$:
$${{T}^{\left( q \right)}}={{P}_{c\left( x \right)}}\times {}_{\left[ 0,q \right]}P,  \quad{{c}_{n}}={{\left( n! \right)}^{-1}}, \quad0\le n<q.$$

E-mail: {evgeniy\symbol{"5F}burlachenko@list.ru}
\end{document}